\documentclass[a4paper,twoside,draft,11pt,righttag]{amsart}
\usepackage{amsxtra}
\usepackage{amsfonts}
\usepackage{amssymb}
\usepackage{pb-diagram}

\theoremstyle{plain}
\newtheorem{prop}{Proposition}[section]
\newtheorem{teo}[prop]{Theorem}

\newtheorem{lema}[prop]{Lemma}

\theoremstyle{definition}

\newtheorem{ejem}[prop]{Example}
\newtheorem{rem}[prop]{Remark}

\theoremstyle{remark}

\numberwithin{equation}{section}

\newcommand{\Z}{\mathbb Z}

\newcommand{\R}{\mathbb R}

\newcommand{\g}{\gamma}
\newcommand{\ld}{\lambda}
\newcommand{\Ld}{\Lambda}
\newcommand{\G}{\Gamma}
\newcommand{\D}{\Delta}
\newcommand{\de}{\delta}

\newcommand{\arr}{\rightarrow}

\newcommand{\sz}{\scriptsize}

\newcommand{\I}{\text{\sl Id}}

\newcommand{\man}{M_\Gamma}

\newcommand{\vep}{\varepsilon}

\newcommand{\on}{\text{O}(n)}
\newcommand{\son}{\text{SO}(n)}
\newcommand{\spin}{\text{Spin}(n)}

\newcommand{\Pinpm}{\text{Pin}^\pm(n)}
\newcommand{\pinpm}{\text{pin}^\pm}
\newcommand{\pinp}{\text{pin}^+}
\newcommand{\pinm}{\text{pin}^-}

\newcommand{\mjh}{M_{j,h}}

\newcommand{\fa}{\mathcal{F}}

\newcommand{\bjh}{B_{j,h}}

\title[spin structures and spectra of $\Z_2^k$-manifolds]{SPIN STRUCTURES
 AND SPECTRA OF $\Z_2^k$-MANIFOLDS.}
\author[R. J. Miatello]{Roberto J. Miatello}
\address{FaMAF--CIEM \\$\!$Universidad Nacional de C\'ordoba\\5000 C\'ordoba, Argentina.}
\email{miatello@mate.uncor.edu, podesta@mate.uncor.edu}
\author[R. A. Podest\'a]{Ricardo A. Podest\'a}
\keywords{flat manifolds, spin structures, isospectrality.}
\thanks{2000 {\it Mathematics Subject Classification.} Primary 58J53, 57R15; \,Secondary 20H15.}
\thanks{Supported by Conicet, Secyt-UNC}
\begin{document}
\bibliographystyle{plain}

%                   ABSTRACT
\begin{abstract}
We give necessary and sufficient conditions for the existence of pin$^\pm$
 and spin structures on Riemannian manifolds with holonomy group $\Z_2^k$.
 For any $n\geq 4$ (resp.\@ $n\geq6$) we give examples of pairs of
 compact manifolds (resp.\@ compact orientable manifolds) $M_1$, $M_2$,
 non homeomorphic to each other, that are Laplace isospectral on
 functions and on $p$-forms for any $p$ and such that $M_1$ admits a
  $\pinpm$ (resp.\@ spin)  structure whereas $M_2$ does not.
\end{abstract}

\maketitle

\section*{Introduction}     \label{s.intro}
Any Riemannian manifold $M$ has naturally associated differential
operators of second order, the Laplacian $\D$ acting on smooth functions
and more generally, the $p$-Laplacian $\Delta_p$ acting on smooth
$p$-forms for $0\leq p\leq n$.
 The Dirac operator $D$ is a first order
operator that can not always be defined. To make this possible, $M$ needs
to have an additional structure: a spin structure, if $M$ is orientable,
and a $\pinpm$ structure, in general. In this case one says that M is spin
or $\pinpm$, respectively.

In this paper we consider a question posed by David Webb, namely, can one
hear the property of being spin on a compact Riemannian manifold? We shall
answer this question in the negative by giving several examples of Laplace
isospectral Riemannian manifolds $M_1,M_2$  such that $M_1$ is spin
(resp.\@ $\pinpm$) but $M_2$ has no spin (resp.\@ $\pinpm$) structure. All
our examples will be isospectral on $p$-forms for $0\le p\le n$
and will be given by {\em $\Z_2^k$-manifolds}, that is, compact Riemannian
manifolds with holonomy group $F\simeq \Z_2^k$. We note that by the
Cartan-Ambrose-Singer theorem, such a manifold is necessarily flat, hence
of the form $M_\G = \G \backslash \R^n$, $\G$ a Bieberbach group.

In one of the main results, Theorem 2.1, we give a parametrization of the
$\pinpm$ or spin structures of $M_\G$, showing that the number is either
$2^r$ for some $r\ge k $ or zero, and  deriving  a simple criterion for
non existence (see Remark 2.3). In Section 3 we apply Theorem 2.1 and this
criterion to construct several isospectral pairs $M,$ $M'$ of
$\Z_2^2$-manifolds of dimensions  $n\ge 4$ (resp.\@ $n \ge 6$), such that
$M$ %one of them
admits a pin$^\pm$ (resp.\@ spin) structure while $M'$ %the other
 does not, thus
giving a negative answer to Webb's question. By increasing dimensions, we
obtain examples of pairs having these same properties and with the extra
condition that both $M,M'$ are K\"ahler (see Remark 3.1).

In the last section we specialize to the case $k=1$, i.e.\@ of
$\Z_2$-manifolds. We show that  any such $M_\G$  has $2^{n-j}$  $\pinpm$
structures for some $0\leq j \leq [\frac{n-1}2]$, with $j$ determined by
the $\Z_2$-action.  If furthermore $M_\G$ is of the so called diagonal
type and orientable, it turns out that $M_\G$ admits $2^n$ spin
structures, as in the case of the $n$-torus (see \cite{Fr}).

\section{Preliminaries}

\subsubsection*{Bieberbach manifolds}
A crystallographic group is a discrete, cocompact subgroup  $\Gamma$ of
the isometry group $I(\R^n)$ of $\R^n$. If $\G$ is torsion-free, then $\G$
is said to be a {\it Bieberbach group}. Such a $\Gamma$ acts properly
discontinuously on $\R^n$, thus $M_\Gamma = \Gamma\backslash\R^n$ is a
compact flat Riemannian manifold with fundamental group $\Gamma$ and
furthermore, any such manifold arises in this way. Since $I(\R^n)\simeq
\on \ltimes \R^n$, any element $\gamma \in I(\R^n)$ decomposes uniquely as
$\gamma = B L_b$, with $B \in \on$ and $b\in \R^n.$ The  translations in
$\Gamma$ form a normal maximal abelian subgroup  of finite index
$L_\Lambda$,  $\Lambda$   a lattice in $\R^n$ which is $B$-stable for
every $BL_b \in \Gamma$. The restriction to $\Gamma$ of the canonical
projection $r:I(\R^n)\rightarrow\text{O}(n)$ given by $BL_b \mapsto B$ is
a homomorphism with kernel $L_{\Lambda}$ and $r(\Gamma)$ is a finite
subgroup of $\on$ isomorphic to $F:= L_\Lambda \backslash \G$. It is
called the {\em holonomy group} of $\Gamma$ and gives the linear holonomy
group of the Riemannian manifold $M_\Gamma$.

A Bieberbach group $\G$ is said to be of {\em diagonal type} (see
\cite{MR2}, Definition 1.3) if there exists an orthonormal $\Z$-basis
$\{\ld_1,\dots,\ld_n\}$ of the lattice $\Lambda$ such that for any element
$BL_b\in\Gamma$, $B\ld_i=\pm \ld_i$ for $1\le i\le n$. These Bieberbach
groups have a rather simple holonomy action, among those with holonomy
group $\Z_2^k$. If $\Gamma$ is of diagonal type, after conjugation of
$\Gamma$ by an isometry,  it may be assumed that  $\Lambda$ is the
canonical lattice and that $b$ lies in $\frac 1{2} \Lambda$ for any
$\gamma=BL_b \in \Gamma$. Thus, any $\gamma \in \Gamma$  can be written
uniquely as $\gamma=BL_{b_o}L_\lambda$, where the coordinates of $b_o$ are
0 or $\frac 12$ and  $\lambda \in \Lambda$ (see \cite{MR2}, Lemma 1.4).

\subsubsection*{Pin and spin groups}
 For a  discussion of the material in this subsection we refer to
 \cite{La}, \cite{Fr2} or \cite{GLP}.
Let $Cl^\pm(n)$ denote the Clifford algebras of $\R^n$ endowed with the
definite quadratic forms $\mp||\cdot||^2$. If $\{e_1, \dots, e_n\}$
denotes the canonical basis of $\R^n$,  then a  basis for $Cl^{\pm}(n)$ is
given by the set $\{ e_{i_1} %e_{i_2}
\dots e_{i_k} \,:\, 1 \leq i_1<\cdots
<i_k \leq n\}$. On $Cl^\pm(n)$ one has the relation  $vw+wv= \pm 2\langle
v,w\rangle$ for any $v,w\in \R^n$, where $\langle \,,\,\rangle$ denotes
the standard inner product. Thus,
\begin{equation} \label{clifford}
\left. \begin{array}{cc}
    e_i e_j  = -e_je_i   \qquad \qquad & \text{ for }  i\ne j
    \text{ for both }  Cl^{\pm}(n), \text{\vspace{1em}}  \medskip \\
    e_i^2   = \pm 1    \qquad  \qquad & \text{ for }
    1\leq i\leq n \text{ in }   Cl^{\pm}(n).
 \end{array} \right.
\end{equation}
 We have compact Lie subgroups, $\text{Pin}^{\pm}(n)$, of the group of units of $Cl^{\pm}(n)$, with $\text{Pin}^{\pm}(n)=\{v_1 \ldots v_{h} \,:\, v_j \in \R^n, \|v_j\| =1 ,\;1 \le j \le h \}$. The connected component of the identity in both cases is isomorphic to
$\text{Spin}(n)= \{v_1 \dots v_{h} \,:\, v_j \in \R^n, \|v_j\| =1 ,\;1 \le j \le h ,\, h \text{  even}\}$,  a compact, simply connected  Lie group for $n\geq3$.

Let $\alpha$ be the canonical involution of $Cl^{\pm}(n)$ given by
$\alpha(v_1\dots v_h)=(-1)^h v_1\dots v_h$. Then,
 we have Lie group epimorphisms
$$\mu_\pm: \text{Pin}^\pm(n) \arr \on$$
 with kernel $\{\pm1\}$, given by $\mu_\pm (v)(x)=\alpha(v)xv^{-1}$ where
  $v\in \text{Pin}^{\pm}(n)$ and $x\in \R^n$.
If $v \in \R^n$, $\|v\| =1$, then $\mu_\pm (v)(x)=-vxv^{-1}= \rho_v(x)$
where $\rho_v$ denotes the orthogonal reflection with respect to the
hyperplane orthogonal to $v$. When restricted to the connected component
of the identity,  $\mu:=\mu_\pm :\spin \simeq \text{Pin}^\pm(n)_o
\rightarrow \text{SO}(n)$ give double coverings.

If $A_j$ is a matrix, for $1\le j \le m$  we will abuse notation by
denoting by $\text{diag}(A_1,\dots, A_m)$ the matrix having $A_j$ in the
``diagonal" position $j$.

Let  $B(t) = \left[\begin{smallmatrix}      \cos t & -\sin t \\
                \sin t & \cos t   \end{smallmatrix}\right]$ with  $t \in \R$ and put
\begin{eqnarray*}
\tau(t_1,\ldots,t_m)=\left\{\begin{array}{ll}\text{diag}(B(t_1),\ldots, B(t_m)), & \quad \text{if }n=2m \\
\text{diag}( B(t_1),\dots, B(t_m),1), & \quad \text{if } n=2m+1.\end{array} \right.
\end{eqnarray*}
 We have that  $T=\left\{\tau(t_1,\ldots,t_m) :t_j \in \R \right\}$ is a maximal torus of $\text{SO}(n)$. A maximal torus of  $\text{Spin}(n)$ is
given by
$$ \tilde{T}=\Big\{ \prod_{j=1}^m \big(\cos t_j + \sin t_j \;e_{2j-1}e_{2j}\big) \;:\; t_j \in \R \Big\}.$$
The restriction $\mu:\tilde{T}\rightarrow T$ is a 2-fold cover and
\begin{equation} \label{e.muinv}
\mu\big(\prod_{j=1}^m (\cos t_j + \sin t_j \;e_{2j-1}e_{2j})\big) = \tau(2t_1,\ldots,2t_m).
\end{equation}

\subsubsection*{Spin structures and $\text{pin}^\pm$ structures.}
If $(M,g)$ is a Riemannian manifold of dimension $n$, let
$\text{B}(M)=\bigcup_{x\in M}\text{B}_{x}(M)$ be the bundle of frames
 on $M$ and $\pi:\text{B}(M)\rightarrow M$ the canonical projection.
 That is, for $x\in M$, $\text{B}_{x}(M)$ is the set of ordered
 orthonormal bases $(v_1,\ldots,v_n)$ of ${T}_x(M)$ and
 $\pi((v_1,\ldots,v_n))=x$.
 $\text{B}(M)$ is a principal $\text{O}(n)$-bundle over $M$
and, if $M$ is orientable, the bundle of oriented frames $\text{B}^{+}(M)$
 is a principal $\text{SO}(n)$-bundle. %over $M$.
A {\em \text{pin}$^\pm$ structure}  on $M$ is a 2-fold cover
$p:\tilde{\text{B}}(M)\rightarrow \text{B}(M)$ that is equivariant and  so
that $\tilde{\pi}:\tilde{\text{B}}(M)\rightarrow M$  is a principal
$\text{Pin}^{\pm}(n)$-bundle
 with $\pi \circ p=\tilde{\pi}$.
Similarly, a  {\it spin structure} on an orientable manifold $M$ is
an equivariant 2-fold cover $p:\tilde{\text{B}}^+(M)\rightarrow \text{B}^+(M)$ where  $\tilde{\pi}:\tilde{\text{B}}^+(M)\rightarrow M$ is a principal $\text{Spin}(n)$-bundle and $\pi \circ p=\tilde{\pi}$.

A manifold in which a spin or a $\text{pin}^\pm$ structure  has been
chosen is called a {\em spin} or a {\em $\text{pin}^\pm$ manifold},
respectively. Note that if $M$ is orientable, any $\pinpm$ structure
on $M$ defines a spin structure and conversely.

We will be interested on spin and $\text{pin}^\pm$ structures on quotients
$\man=\Gamma\backslash \R^n$, where $\Gamma$ is a Bieberbach group. If
$M=\R^n$, we have that $\text{B}(\R^n)=\R^n\times\text{O}(n)$, thus
clearly $\R^n\times \text{Pin}^\pm(n)$ are principal
$\text{Pin}^\pm(n)$-bundles and the maps   $\I \times \mu_\pm:\R^n\times
\text{Pin}^\pm(n) \rightarrow \R^n\times \text{O}(n)$  are equivariant
2-fold covering maps.
 Similarly, we have that
$\R^n\times \spin$ is a principal $\text{Spin}(n)$-bundle and an
equivariant $\text{2-fold}$ cover of $\text{B}^+(\R^n)=\R^n\times
\text{SO}(n)$.
 Thus we have  spin and \text{pin}$^\pm$ structures on $\R^n$ and since
 $\R^n$ is contractible these are the only such structures.
Now, if $\Gamma$ is a Bieberbach group we have a left action of $\Gamma$
on $\text{B}(M)$ given by $\gamma\cdot(x,(w_1,\ldots,w_n)) = (\gamma
x,(\gamma_\ast w_1,\ldots,\gamma_\ast w_n))$. If $\gamma=B L_{b}$ then
$\gamma_\ast w_j=w_j B$. Fix $(v_1,\dots,v_n) \in \text{B}(M)$. Since
$(w_1,\ldots,w_n)=(v_1k,\ldots,v_nk)$ for some $k\in \text{O}(n)$, we see
that $\gamma_\ast w_j=(v_jk)B= v_j(B k)$, thus the action of $\Gamma$ on
$\text{B}(M)$  corresponds to the action of $\Gamma$ on $\R^n\times
\text{O}(n)$ given by $\gamma\cdot(x,k)=(\gamma x,B k)$.

Now assume that there is a group homomorphism $\varepsilon:\Gamma
\rightarrow \text{Spin}(n)$ (resp.\@ $\varepsilon_\pm:\Gamma \rightarrow
\Pinpm$) such that $\mu(\varepsilon (\gamma))=r(\gamma)$ (resp.\@
$\mu_\pm(\varepsilon_\pm(\gamma))=r(\gamma)$). In this case we can lift
the left action of $\Gamma$ on $\text{B}^+(\R^n)$
 (resp.\@ on $\text{B}(\R^n)$)  to
  $\tilde{\text{B}}^+(\R^n)=\R^n\times \spin$
(resp.\@ to $\tilde{\text{B}}(\R^n)=\R^n\times \Pinpm$) via
 $\gamma\cdot(x,\tilde k)=(\gamma x,\vep(\gamma) \tilde k)$.
 Thus we have the spin structure
$$\begin{diagram} \node{\Gamma\backslash (\R^n\times \spin)} \arrow[2]{e,t}{\overline{\I\times \mu}} \arrow{se} \node{} \node{\Gamma \backslash (\R^n\times \son)} \arrow{sw}
 \\ \node{} \node{\Gamma\backslash \R^n} \node{}
\end{diagram}$$
for $M_\G$ since $\Gamma \backslash \text{B}(\R^n)=\text{B}(\G \backslash \R^n)$  and 
$\overline{\I\times \mu}$ is equivariant. Similarly for the $\pinpm$ structures.

In this way, for each homomorphism $\varepsilon$ or $\vep_\pm$ as above, we obtain   a spin or a  pin$^\pm$ structure on $M_\Gamma$, respectively. It turns out that all spin  and \text{pin}$^\pm$ structures on $M_\Gamma$ are obtained in this manner (see \cite{Fr2}, \cite{La}).

The $n$-torus admits $2^n$  spin structures. Indeed,  if $T_\Lambda
=\Lambda\backslash \R^n$, and $\ld_1, \dots, \ld_n$ is a $\Z$-basis of
$\Lambda$, then a homomorphism $\vep$ as above is determined by the
$n$-tuple $\vep (L_{\ld_i})=\delta_i \in \{\pm1\}$, for $1\le i \le n$
(see \cite{Fr}).
We shall show in Section 4 that  this is still the number of such
 structures for flat manifolds with holonomy group $\Z_2$ which are of
 diagonal type.

\section{Spin and $\pinpm$ structures on $\Z_2^k$-manifolds.}
 In this section we study  the  existence of  $\text{pin}^\pm$ structures
  on $\Z_2^k$-manifolds, showing  that the number of such structures is
  either 0 or $2^r$ for some  $r\geq k$.  
As an application, in the next section we will construct many examples
 of $\Z_2^2$-manifolds for any $n\ge 4$,  having $\pinp$ structures but
 no $\pinm$ structures (and conversely) or else,  having neither of them.

Let $\G$ be a Bieberbach group with holonomy group $F\simeq \Z_2^k$, $1\le
k \le n-1$, and translation lattice $\Ld$. Then  $M_\G = \G\backslash
\R^n$ with $\G=\langle \g_1,\dots,\g_k, \Ld \rangle$ where
$\g_i=B_iL_{b_i}$, $B_i \in \on$, $b_i \in \R^n$,   $B_i\Ld=\Ld$,
$B_i^2=\I$ and $B_i B_j=B_j B_i$, for each $1\leq i,j\leq k$.

Assume there is a $\pinpm$ structure on $M_\G$, that is, a group
homomorphism $\vep_\pm: \G \rightarrow \Pinpm$ such that $\mu_\pm\circ\vep_\pm=r$. Then, necessarily $\vep_\pm(L_{\ld})\in \{\pm1\}$, for $\ld \in \Ld$.
Thus, if $\ld_1,\dots,\ld_n$ is a $\Z$-basis of $\Ld$ and we set  $\delta_i:= \vep_\pm(L_{\ld_i})$, for every $\ld =\sum_im_i \ld_i \in \Ld$ with $m_i \in \Z$, we have $\vep_\pm(L_\ld) = \prod_i \de_i^{m_i}=\prod_{m_i \text{odd}} \delta_i$.  

If $\g=BL_b\in \G$ we will fix  a distinguished (though arbitrary) element in  $\mu^{-1}_\pm(B)$, denoted by $u_\pm(B)$.
If $M_\G$ is orientable, we write  $u(B):=u_\pm(B)$.
Thus, if $\g =BL_{b}\in \G$, then
\begin{equation}\label{sigma}
\vep_\pm(\g) = \sigma \, u_\pm(B),
\end{equation}
 where $\sigma \in \{\pm 1\}$ depends on $\g$ and on the choice of $u_\pm(B)$.

Let $\G=\langle \g_1,\dots,\g_k,\Ld \rangle$.
The morphism $\vep_\pm$ is determined by its action on the generators of $\G$, that is, by the $(n+k)$-tuple
\begin{equation}\label{tuple}
\begin{split}
(\delta_1,\ldots,\delta_n, \sigma_1 u_\pm(B_1),\ldots,\sigma_k u_\pm(B_k))\\ \text{ or } \qquad
(\delta_1,\dots,\delta_n,\sigma_1,\dots,\sigma_k)\in \{\pm 1\}^{n+k}
\end{split}
\end{equation}
where $\delta_i = \vep_\pm(L_{\ld_i})$ and $\sigma_i$ is defined by the equation $\vep_\pm(\g_i)=\sigma_i u_\pm(B_i)$, for $1 \le i \le k$.

Now, since $\vep_\pm$ is a homomorphism, for any $\g =BL_b \in \G, \ld \in
\Ld$ we  have
\begin{equation*}
\vep_\pm(L_{B\ld})=\vep_\pm(\g L_\ld \g^{-1})=\vep_\pm(\g) \vep_\pm(L_\ld)
\vep_\pm(\g^{-1}) =\vep_\pm(L_\ld).
\end{equation*}
Therefore we see that if $\vep_\pm$ is a $\pinpm$ structure on $M_\G$,
since $\g^2 \in L_{\Ld}$, then the character ${\vep_\pm}_{|\Ld}$ must
satisfy the following conditions for any $\g=BL_b \in \G$:
\begin{equation}\label{conditions}
\left.\begin{array}{ll} (\vep_1)& \qquad \vep_\pm (\g^2)= \vep_\pm (\g)^2
= u_\pm^2(B) \medskip \\ (\vep_2) &\qquad \vep_\pm(L_{(B-\I)\ld})=1,\quad
\ld \in \Ld.
\end{array}\right.
\end{equation}
We thus set
\begin{equation}\label{veplambda}
\hat{\Ld} (\G):=\{ \chi \in \text{Hom}(\Ld,\{\pm1\}) \,:\, \chi
\text{ satisfies } (\vep_1) \text{ and } (\vep_2) \}.
\end{equation}

The next result gives a parametrization of  the  $\pinpm$ structures
 $\vep_\pm$ for $\man$.

\begin{teo} \label{z2spinst}
 If $\G=\langle \g_1,\ldots,\g_k, \Ld \rangle$ is a Bieberbach group with
  holonomy group $\Z_2^k$ and $\sigma_1, \dots, \sigma_k$ are as in
  (\ref{tuple}), then the map $\vep_\pm \!\mapsto \!
  ({\vep_\pm}_{|\Ld},\sigma_1,\ldots, \sigma_k)$   defines a bijective
   correspondence between the $\pinpm$ structures on $\man$ and the set
    $\hat{\Ld} (\G)\times \{\pm1\}^k$.
The number of pin$^\pm$  structures on $\man$ is either $0$ or $2^r$ for
 some $r \ge k$.
\end{teo}

\begin{proof} %In the course of the proof
We shall write $\vep, \mu, u(B)$ in place of $\vep_\pm$, $\mu_\pm$,
$u_\pm(B)$, for simplicity.

Any element $\g \in \G$ can be written as a product of generators $\g_i
=B_iL_{b_i}$ and $\ld \in \Ld$.  After reordering, by normality of $\Ld$ in $\G$ and since
$B_i^2= \I$, we see that $\g$ can be written uniquely as
\begin{equation}\label{unique}
\g = \g_{i_1} \ldots \g_{i_r} L_\ld, \quad \text{ with }1 \le i_1 <\dots < i_r\le k,\, \ld \in \Ld.
\end{equation}

Given $ \vep \in \hat{\Ld}(\G)$  and for any choices of $\vep(\g_i) \in
\mu^{-1}(B_i)$, $1\le i \le k$,  we define (in the notation of
(\ref{unique})) for  $\g \in \G$:
\begin{equation}\label{defvarep}
\vep(\g)= \vep(\g_{i_1}) \ldots \vep(\g_{i_r}) \vep(L_\ld).
\end{equation}
Thus, we get a well defined map $\vep : \G \rightarrow \Pinpm$ such that
$\mu \circ \vep = r$ and we claim it is a homomorphism. For this purpose
we need to show that
\begin{equation}\label{multip}
\vep(\g_{i_1} \dots \g_{i_r}L_\ld \g_{j_1} \dots \g_{j_t} L_{\ld'})
= \vep(\g_{i_1} \dots \g_{i_r}L_\ld) \vep(\g_{j_1} \dots \g_{j_t} L_{\ld'}),
\end{equation}
for any $i_1<\dots < i_r, \, j_1<\dots<j_t$ and $\ld, \ld' \in \Ld$.

We first note that we may leave out $\ld, \ld'$ in (\ref{multip}). Indeed,
assume that for  $\g,\g'\in \G$ one has $\vep(\g \g')= \vep(\g)
\vep(\g')$. Then, by $(\vep_2)$
\begin{eqnarray*}
\vep(\g L_\ld \g'L_{\ld'}) & = & \vep(\g \g'L_{B\ld+\ld'})
 = \vep(\g \g')\vep(L_{B\ld+\ld'}) \\
& = & \vep(\g) \vep(\g') \vep(L_{\ld}) \vep(L_{\ld'}) =  \vep(\g L_\ld)
\vep(\g'L_{\ld'}).
\end{eqnarray*}

As a step in the proof of (\ref{multip}) (with $\ld=\ld'= 0$) we will first show that
\begin{equation}\label{homotwo}
\vep(\g_i \g_j)=\vep(\g_i)\vep(\g_j), \text{ for any } i,j.
\end{equation}
This follows from the definition of $\vep$, if $i<j$, and from condition
($\vep_1$), if $i=j$. We thus assume that $j<i$.
Then we may write $\g_i \g_j = \g_j \g_i [{\g_i}^{-1}, {\g_j}^{-1}]$.
Since $[{\g_i}^{-1}, {\g_j}^{-1}] \in \Ld$, by the definition of $\vep$
\begin{equation}\label{vepji}
\vep(\g_i \g_j) = \vep(\g_j) \vep(\g_i) \vep( [{\g_i}^{-1}, {\g_j}^{-1}]).
\end{equation}
Note that (\ref{vepji}) will equal $ \vep(\g_i) \vep(\g_j)$ if and only if
it holds the relation
\begin{equation} \label{vepcommut}
\vep( [{\g_i}^{-1}, {\g_j}^{-1}]) = [\vep({\g_i}^{-1}),
\vep({\g_j}^{-1})].
\end{equation}

To show (\ref{vepcommut}), we have by condition ($\vep_1$) that
\begin {equation}\label{vepijsquare}
\vep ((\g_j \g_i)^2)= \vep (\g_j \g_i)^2= \vep (\g_j) \vep(\g_i)\vep (\g_j) \vep(\g_i).
\end{equation}
On the other hand
\begin{equation}\label{vepijsquare2}
\begin{split}
\vep ((\g_j \g_i)^2)%&=& \vep (\g_j  \g_i \g_j \g_i)\\
& \;=\; \vep (\g_j \g_j \g_i [{\g_i}^{-1}, {\g_j}^{-1}] \g_i) \\
& \;=\; \vep (\g_j^2 \g_i^2 ({\g_i}^{-1}[{\g_i}^{-1}, {\g_j}^{-1}] \g_i))\\
& \;=\; \vep({\g_j}^2) \vep({\g_i}^2) \vep({\g_i}^{-1}[{\g_i}^{-1}, {\g_j}^{-1}] \g_i)\\
& \;=\; \vep({\g_j})^2 \vep({\g_i})^2 \vep([{\g_i}^{-1}, {\g_j}^{-1}]).
\end{split}
\end{equation}
In the last  equality we have used condition ($\vep_2$) and the fact that
commutators lie in $\Lambda$.

Now, by combining (\ref{vepijsquare}) and (\ref{vepijsquare2}) we obtain
(\ref{vepcommut}), hence (\ref{homotwo}) follows.

\smallskip
In the general case, (\ref{multip}) can be proved by an inductive
argument.

     Let first $t=1$, $r$ arbitrary. The case $r=1$ is (\ref{homotwo}), so
assume $r>1$. If $j_1>i_r$, then the assertion is clear by the definition
of $\vep$, while if $j_1= i_r$, we may use ($\vep_1$) and induction. We
thus assume that there is $\alpha$ such that $i_{\alpha-1} \le j_1 <
i_\alpha$. Actually, we shall take %assume that
$i_{\alpha-1} < j_1 $. The
proof
when $i_{\alpha-1} = j_1$ is similar, but simpler.

If we set $u=\g_{i_\alpha}\cdots \g_{i_r}$ then
\begin{eqnarray*}
\vep(\g_{i_1}\cdots \g_{i_r} \g_{j_1}) &=& \vep(\g_{i_1}\cdots \g_{i_{\alpha-1}} \g_{j_1} 
u [u^{-1}, \g_{j_1}^{-1}])
\\
&=& \vep(\g_{i_1})\cdots \vep(\g_{i_{\alpha-1}}) \vep( \g_{j_1}) \vep(u)  \vep([u^{-1}, \g_{j_1}^{-1}]) \\
(\text{by }(\ref{vepcommut}))&=& \vep(\g_{i_1})\cdots \vep(\g_{i_{\alpha-1}})  \vep(u)  \vep( \g_{j_1}) \\
&=& \vep(\g_{i_1}\cdots \g_{i_r}) \vep(\g_{j_1}).
\end{eqnarray*}
The argument  for arbitrary $t$ is quite  similar and will be omitted.
\end{proof}

\begin{rem} \label{cond.b}
For manifolds of diagonal type, condition ($\vep_2$)   always holds, since $(B-\I)\Ld \subset 2\Ld$ for any $BL_b \in \G$.   More generally, for manifolds whose holonomy  representation decomposes as a sum of  integral representations of rank $\leq 2$, condition ($\vep_2$)  can be expressed in simple terms.

In Section 4 we will study in more detail the case of $\Z_2$-manifolds, showing in particular that $\pinpm$ structures can always be defined in this case.
\end{rem}

\begin{rem} \label{criterion} The previous theorem shows that there are restrictions for a
$\Z_2^k$-manifold $\man$ to carry a $\pinpm$ structure. As a consequence, one has the following simple criterion:

{\it Suppose there exist $\g=BL_b, \g'=B'L_{b'} \in \G$ with $\g^2={\g'}^2$ and such that for $u_+(B) \in \mu_+^{-1}(B)$ and $u_+(B') \in \mu_+^{-1}(B')$ one has $u_+(B)^2=-u_+(B')^2$.
Then $\man$ can not admit a $\pinp$ structure.}

Indeed, such a structure $\vep_+$ would have to satisfy $\vep_+(\g)=\pm u_+(B)$, $\vep_+(\g')=\pm u_+(B')$ and $\vep_+(\g^2)=\vep_+({\g'}^2)$, that is, $u_+(B)^2=u_+(B')^2$ against our assumption.

The same criterion, with the obvious changes, is valid for non existence of $\pinm$ structures,
or spin structures in the orientable case.
\end{rem}

\begin{rem}\label{doubling}
In contrast with Remark \ref{criterion},  by applying the doubling
procedure in \cite{DM2}, we may obtain spin Bieberbach manifolds of
diagonal type with holonomy group $\Z_2^k$, for any $k\ge 1$. Indeed, let
$\G=\langle \g_1,\ldots,\g_k, L_{\Ld} \rangle$ be an $n$-dimensional
Bieberbach group of diagonal type with holonomy group $\Z_2^k$. Define
$\text{d}\G := \langle d\g_1,\ldots,d\g_k, L_{\Ld \oplus \Ld} \rangle $
where $d\g := \left[ \begin{smallmatrix} B & 0 \\ 0 & B \end{smallmatrix}
\right]  L_{(b,b)}$ if $\g=BL_b\in \G$ (see Definition 3.1 in \cite{DM2}).
Thus, $\text{d}\G$ is a Bieberbach group of dimension $2n$ with holonomy
group $\Z_2^k$. The manifold $M_{\text{d}\G}=\text{d}\G\backslash \R^{2n}$
is an orientable K\"ahler flat manifold of diagonal type. If we apply this
procedure twice, then the manifold $M_{{\text{d}^2}\G}$ is hyperk\"ahler
(see Proposition 3.2 in \cite{DM2}). It turns out that this
$4n$-dimensional manifold is always spin. Indeed, in the notation of Lemma
\ref{preimages} in the next section, since $h\in 4\Z$ for
${\text{d}^2}\G$,  we have that $u^2(B)= u^2_{0,h}=1$ by (\ref{squares}).
Hence, condition ($\vep_1$) takes the form $\vep(\g^2)=1$ for any $\g \in
\G$. Therefore, spin structures can always be defined for
$M_{{\text{d}^2}\G}$, for example we may take any of the $2^{k}$
homomorphisms $\vep :\G \rightarrow \spin $ such that $\vep_{|\Ld} \equiv
1$.
\end{rem}

\section{Spin structures on some  isospectral pairs.}
In this section we will construct several isospectral pairs $\{M,M'\}$ of
 $\Z_2^2$-manifolds of dimension 4
by using the results in
\cite{MR2}, and we will determine the $\pinpm$ or spin structures, showing
that, for some of them, $M$ has a $\pinpm$ or a spin structure, while $M'$
does not. The main result is given in Theorem \ref{main}. In the proof, we
will need to know some preimages in $\Pinpm$ by $\mu_\pm$, as well as
their squares.

Set $J:=\left[ \begin{smallmatrix} 0 & 1 \\ 1 & 0 \end{smallmatrix} \right]$. For each $0\leq j,h < n$, we set
\begin{equation}\label{Bjh}
B_{j,h}=\text{diag}(\underbrace{J,\dots,J}_j,\underbrace{-1,\dots,-1}_h,\underbrace{1,\dots,1}_l),
\end{equation}
where $n=2j+h+l$, $j+h\not=0$ and $l\geq1$.

\begin{lema} \label{preimages}   Let $B_{j,h}$ be as in (\ref{Bjh}) and  let $\mu_\pm:\Pinpm \rightarrow \on$ be the canonical covering maps. If we set
\begin{equation}\label{ujh}
u^\pm_{j,h}:=(\tfrac{\sqrt 2}2)^j (e_1 -e_2) \ldots (e_{2j-1}-e_{2j})\, e_{2j+1} \cdots e_{2j+h},
\end{equation}
then
$\mu_+^{-1}(B_{j,h})= \{\pm u^+_{j,h}\}$, $\mu_-^{-1}(B_{j,h})= \{\pm u^-_{j,h}\}$ and furthermore
\begin{equation}
\begin{split}\label{squares}
(u_{j,h}^+)^2 & = (-1)^{jh}(-1)^{[\frac{j}2]}(-1)^{[\frac{h}2]} \\
(u_{j,h}^-)^2 & =   (-1)^{jh}(-1)^{[\frac{j+1}2]}(-1)^{[\frac{h+1}2]}.
\end{split}
\end{equation}
In particular,  $(u_{0,h}^+)^2 = (-1)^{[\frac{h}2]}$
and $(u_{0,h}^-)^2 = (-1)^{[\frac{h+1}2]}$.
 If $B_{j,h} \in SO(n)$, i.e.\@ if $j+h$ is even, then
$u_{j,h}^2 = (-1)^{\frac{j+h}2}$.

If $B \in \on$ is conjugate to $\bjh$, and $u_\pm(B)\in \mu_\pm^{-1}(B)$, then $u^2_\pm(B)=(u_{j,h}^\pm)^2$.
\end{lema}

\begin{proof} Since $\mu_\pm (e_i)=\rho_{e_i}=\text{diag}(1,\dots,1,-1 ,1,\dots,1)$ with $-1$  in the $i$-th position,
it is clear that  $\mu_+^{-1}(B_{0,h})=\mu_-^{-1}(B_{0,h})=\big\{\pm e_1 \ldots e_h \big\}$.
If $n=2$, 
we may write $J$  as a product
$
J=\left[ \begin{smallmatrix} -1&0\\0&1 \end{smallmatrix} \right] \left[ \begin{smallmatrix} 0&-1\\1&0 \end{smallmatrix} \right]$.
Hence, using (\ref{e.muinv}), and $\mu_\pm (e_i) = \rho_{e_i}$, we get that $\mu_+ ^{-1}(J)=\mu_- ^{-1}(J) =\big \{ \pm e_1(\cos(\tfrac{\pi}4)+\sin(\tfrac{\pi}4)e_1 e_2) \big \}=\big\{ \pm\tfrac{\sqrt2}2 (e_1 - e_2) \big\}$. Arguing similarly for arbitrary $n$, the first assertion in the lemma follows.

On the other hand one computes, using (\ref{clifford}), that both $(e_1\dots e_h)^2$ and  $2^{-h}((e_1 -e_2) \cdots (e_{2h-1}-e_{2h}))^2$ equal $(-1)^{[\frac{h}2]}$ in $Cl^+(n)$ and $(-1)^{[\frac{h+1}2]}$ in $Cl^-(n)$, respectively. This implies equations  (\ref{squares}).

Now, suppose $B=C\bjh C^{-1}$ with $C\in\on$. If $u_+(C)
 \in \mu_+(C)^{-1}$, then  $u_+(B)=\pm u_+(C) u_{j,h}^+ u_+(C)^{-1}$ and hence  $u^2_+(B)=u_+(C) (u_{j,h}^+)^2 u_+(C)^{-1} = (u_{j,h}^+)^2$. The verification for $u^2_-(B)$ is identical.
\end{proof}

\medskip

We now consider  some pairs of 4-dimensional $\Z_2^2$-manifolds
 $\{M_i,M_i'\}$, $1\le i \le 5$,
where $M_i=\G_i\backslash \R^4$, $M_i'=\G_i'\backslash \R^4$  and the
groups $\G_i=\langle \g_1, \g_2, \Ld  \rangle$, $\G_i'=\langle \g_1',
\g_2', \Ld  \rangle$ are given in Table 1, where $\g_i=B_i L_{b_i}$,
$\g_i'=B_i L_{b_i'},$  $i=1,2$, $B_3=B_1B_2$, $b_3=B_2 b_1 + b_2$,
$b_3'=B_2' b_1' + b_2'$ and $\Ld=\Z e_1 \oplus \ldots \oplus \Z e_n$ is
the canonical lattice. Furthermore, we take  $B_i = B'_i$. In all cases
the matrices $B_i$ are diagonal and are written as column vectors. We
indicate
 the translation vectors $b_i, b'_i$ also as column vectors, leaving out
  the coordinates that are equal to zero. We will also use  the pair
$\{\tilde{M_1},\tilde{M_1'}\}$ of $\Z_2^2$-manifolds of dimension 6
obtained from the pair $\{M_1,M_1'\}$  by adjoining the characters
$(-1,1,-1)$ and $(1,-1,-1)$ to $B_i$, $1 \le i \le 3$, and keeping $b_i$, $b'_i$
unchanged.

\medskip
\begin{center}
{\sc Table 1}

\medskip
$\left.\begin{array}{c} \{M_1,M_1'\}\\
\{\tilde{M}_1,\tilde{M}'_1\}
\end{array}\right.$
 \qquad
\begin{tabular}{|rcc|rcc|rcc|}  \hline    $B_1$ &  $L_{b_1}$  & $L_{b_1'}$    & $B_2$ &  $L_{b_2}$ & $L_{b_2'}$  & $B_3$ &  $L_{b_3}$ & $L_{b_3'}$   \\ \hline
 1 && & 1 &$\text{{\scriptsize 1/2}}$ & $\text{{\scriptsize 1/2}}$ & 1 &$\text{{\scriptsize 1/2}}$ &$\text{{\scriptsize 1/2}}$ \\
      1 && $\text{{\scriptsize 1/2}}$ & 1 & $\text{{\scriptsize 1/2}}$ & & 1 &$\text{{\scriptsize 1/2}}$ &$\text{{\scriptsize 1/2}}$  \\
 1 &  $\text{{\scriptsize 1/2}}$ & &  -1 && & -1 &  $\text{{\scriptsize 1/2}}$ &  \\
 -1 &  & & 1 & &$\text{{\scriptsize 1/2}}$ & -1 & &$\text{{\scriptsize 1/2}}$  \\
 \hline
-1 &  & &  1 & & & -1 &  &  \\ 1 &  & & -1 & & & -1 & & \\  \hline
\end{tabular}
\end{center}

\medskip

\begin{center}
$\left. \begin{array}{c} \{M_2,M_2'\} \end{array} \right. \qquad $
\begin{tabular}{|rcc|rcc|rcc|}  \hline    $B_1$ &  $L_{b_1}$  & $L_{b_1'}$    & $B_2$ &  $L_{b_2}$ & $L_{b_2'}$  & $B_3$ &  $L_{b_3}$ & $L_{b_3'}$   \\ \hline
1 && & 1 & &$\text{{\scriptsize 1/2}}$ & 1 & &$\text{{\scriptsize 1/2}}$ \\
1 && $\text{{\scriptsize 1/2}}$ & 1 & $\text{{\scriptsize 1/2}}$ & $\text{{\scriptsize 1/2}}$& 1 &$\text{{\scriptsize 1/2}}$ &  \\
1 &  $\text{{\scriptsize 1/2}}$ & &  -1 && & -1 &  $\text{{\scriptsize 1/2}}$ &  \\
-1 &  & & 1  &$\text{{\scriptsize 1/2}}$ & & -1  &$\text{{\scriptsize 1/2}}$ & \\
\hline
\end{tabular}
\end{center}

\medskip

\begin{center}
 $ \left. \begin{array}{c} \{M_3,M_3'\} \end{array} \right.$ \qquad
\begin{tabular}{|rcc|rcc|rcc|}  \hline  $B_1$ &  $L_{b_1}$  & $L_{b_1'}$    & $B_2$ &  $L_{b_2}$ & $L_{b_2'}$  & $B_3$ &  $L_{b_3}$ & $L_{b_3'}$   \\ \hline
1 & & &           -1 & &   & -1 & & \\
1 & & $\text{{\scriptsize 1/2}}$ & -1 &   & & -1 & & $\text{{\scriptsize 1/2}}$  \\
-1 &  & &      -1 & $\text{{\scriptsize 1/2}}$ &       &  1 &  $\text{{\scriptsize 1/2}}$ &  \\
1 & $\text{{\scriptsize 1/2}}$  & & 1 & $\text{{\scriptsize 1/2}}$ & $\text{{\scriptsize 1/2}}$ &          1 & & $\text{{\scriptsize 1/2}}$  \\
\hline
\end{tabular}
\end{center}

\medskip

\begin{center}
 $ \left. \begin{array}{c} \{M_4,M_4'\} \end{array} \right.$ \qquad
\begin{tabular}{|rcc|rcc|rcc|}  \hline  $B_1$ &  $L_{b_1}$  & $L_{b_1'}$    & $B_2$ &  $L_{b_2}$ & $L_{b_2'}$  & $B_3$ &  $L_{b_3}$ & $L_{b_3'}$   \\ \hline
1 &$\text{{\scriptsize 1/2}}$  & &           -1 & &   & -1 &$\text{{\scriptsize 1/2}}$  & \\
1 &$\text{{\scriptsize 1/2}}$  & $\text{{\scriptsize 1/2}}$ & -1 &   & & -1 &$\text{{\scriptsize 1/2}}$  & $\text{{\scriptsize 1/2}}$  \\
-1 &  & &      -1 & & $\text{{\scriptsize 1/2}}$  &  1 & & $\text{{\scriptsize 1/2}}$   \\
1 & & $\text{{\scriptsize 1/2}}$ & 1 & $\text{{\scriptsize 1/2}}$ & $\text{{\scriptsize 1/2}}$ &          1 & $\text{{\scriptsize 1/2}}$ & \\
\hline
\end{tabular}
\end{center}

\medskip
\begin{center}
 $ \left. \begin{array}{c} \{M_5,M_5'\} \end{array} \right.$ \qquad
\begin{tabular}{|rcc|rcc|rcc|}  \hline  $B_1$ &  $L_{b_1}$  & $L_{b_1'}$    & $B_2$ &  $L_{b_2}$ & $L_{b_2'}$  & $B_3$ &  $L_{b_3}$ & $L_{b_3'}$   \\ \hline
-1 & & &      1 & &  $\text{{\scriptsize 1/2}}$ & -1 & &$\text{{\scriptsize 1/2}}$ \\
-1 & & &   -1 & $\text{{\scriptsize 1/2}}$ & $\text{{\scriptsize 1/2}}$& 1 &$\text{{\scriptsize 1/2}}$ & $\text{{\scriptsize 1/2}}$    \\
1  & & $\text{{\scriptsize 1/2}}$  &      -1 && & -1 & &  $\text{{\scriptsize 1/2}}$   \\
1  & $\text{{\scriptsize 1/2}}$ &  & 1  & $\text{{\scriptsize 1/2}}$ & & 1  & & \\
\hline
\end{tabular}
\end{center}

\medskip
We observe that only $M_5,M_5',\tilde{M}_1,{\tilde{M}_1}'$ are
orientable.

\smallskip

In order to show the isospectrality of these pairs we will need to recall some known results.

For $BL_b \in \G$ set $n_B:= \text{dim}(\R^n)^B=|\{ 1\le i \le n: Be_i=e_i \}|$
and
\begin{equation}
 n_B(\tfrac 12) := |\{1 \le i \le n: Be_i=e_i \text { and } b\cdot e_i= \tfrac 12 \}|.
\end{equation}
If $0\le t \le d \le n$, the {\em Sunada numbers} for $\G$  are defined
by
\begin{equation}\label{sunada}
c_{d,t}(\G) :=\big|\big\{BL_b \in \G : n_B=d \text{ and } n_B(\tfrac 12) =t \big \}\big|.
\end{equation}
In \cite{MR2}, Theorem 3.3, it is shown that the equality of the Sunada
numbers $c_{d,t}(\G)=c_{d,t}(\G')$  for every $d,t$,  is equivalent to the
validity of the conditions in Sunada's theorem (see \cite{Su}) for $M_\G$
and $M_{\G'}$. In particular this implies that $M_\G$ and $M_{\G'}$ are
isospectral on $p$-forms for $0\le p\le n$. This method was used in
\cite{MR} and \cite{MR3} to prove the isospectrality of the pairs
$M_5,M_5'$ and $M_2,M'_2$ respectively. Also, the method of adding
characters and keeping isospectrality was also used in \cite{MR}.

We are now in a position to state the main result in this paper.

\begin{teo}\label{main}
The pairs $M_i,M_i'$, $1\le i \le 5$, and $\tilde{M}_1, {\tilde{M}_1}'$
are pairwise isospectral.

The number of $\pinpm$ and spin structures on $M_i,M_i'$, $1\le i \le 5$,
and $\tilde{M}_1, {\tilde{M}_1}'$ are  given in the following table.
\renewcommand{\arraystretch}{1.2}
\begin{table}[h]
\label{numberofstr}
\begin{center}
\begin{tabular}{||c||c|c||c|c||c|c||c|c||c|c||c|c||}
\hline
Pairs & $M_1$ & $M_1'$ & $\tilde{M}_1$ & $\tilde{M}'_1$ & $M_2$ & $M_2'$ & $M_3$ & $M_3'$ & $M_4$ & $M_4'$ & $M_5$ & $M_5'$  \\ \hline
$\pinp$ & -- & $2^3$ & -- & $2^5$ & $2^3$ & -- & -- & $2^4$ & $2^4$ & $2^3$ & $2^4$ & $2^3$ \\ \hline
 $\pinm$ & $2^4$ & $2^3$ &-- & $2^5$ & $2^3$ & -- & -- & -- & -- & $2^3$ & $2^4$ & $2^3$ \\ \hline
 spin & -- & -- & --& $2^5$ & -- & -- & -- & -- & -- & -- & $2^4$ & $2^3$ \\ \hline
\end{tabular}
\end{center}
\end{table}

The various isospectral pairs in the table show  that one can not hear the
existence of $\pinpm$ or spin structures on a compact Riemannian manifold.
\end{teo}

\begin{proof}
 Since all manifolds are of diagonal type, to show that these pairs are
isospectral it suffices to check the equality of the Sunada numbers (see
(\ref{sunada})). It is easy to see from Table 1 that the non trivial
Sunada numbers, besides $c_{4,0}=1$ corresponding to the identity, are:
 $c_{2,2}=c_{3,1}=c_{3,2}=1$ for $M_1$ and $M_1'$;
$c_{2,2}=c_{4,1}=c_{4,2}=1$ for $\tilde{M}_1$ and $\tilde{M}_1'$;
$c_{2,1}=c_{3,1}=c_{3,2}=1$ for $M_2$ and $M_2'$;
$c_{1,1}=c_{2,1}=c_{3,1}=1$ for $M_3$ and $M_3'$;
$c_{1,1}=c_{2,1}=c_{3,2}=1$ for $M_4$ and $M_4'$;
 and $c_{2,1}=3$ for $M_5$ and $M_5'$. Thus, it follows that all pairs
 $M_i,{M_i}'$, $1\le i\le 5$, and $\tilde{M}_1, \tilde{M}_1'$ are
 isospectral on functions.
\medskip

 We shall now use Theorem \ref{z2spinst} to determine the spin and $\pinpm$ structures on $M_1,M'_1,\ldots, M_5, M'_5,  \tilde{M}_1$ and $\tilde{M}'_1$. By Remark \ref{cond.b} we need only look at condition ($\vep_1$).

We first look at the pair $M_1, M'_1$. We have that
$$\g_1^2=L_{e_3}, \: \g_2^2=L_{e_1+e_2}=\g_3^2; \qquad {\g_1'}^2=L_{e_2}, \: {\g_2'}^2=L_{e_1+e_4}, \: {\g_3'}^2=L_{e_1+e_2}.$$
By (\ref{e.muinv}) and Lemma \ref{preimages}:
\begin{gather*} u_\pm^2(B_1)=u_\pm^2(B_1')=(\sigma_1 e_4)^2= \pm1, \qquad u_\pm^2(B_2)=u_\pm^2(B_2')= (\sigma_2 e_3)^2=\pm1, \\ u_\pm^2(B_3)=u_\pm^2(B_3')= (\sigma_3 e_3e_4)^2=-1
\end{gather*}
with $\sigma_i \in \{\pm 1\}$.
By the criterion in Remark \ref{criterion}, it follows that  $M_1$ has no $\pinp$ structures, since $\g_2^2= \g_3^2$ and $u_+^2(B_2)= 1$ while  $u_+^2(B_3)=-1$.

Furthermore, by the previous equations, if  $\delta_i =\vep_\pm(L_{e_i})$, condition ($\vep_1$) gives $\delta_3=\pm1$, $\delta_1\delta_2=\pm1$ and $\delta_1\delta_2=-1$. The  last two equations are not compatible for $Cl^+(n)$, hence we see again that $M_1$ does not admit $\pinp$ structures. However, it has $2^4 $ $\pinm$ structures given by
$$\vep_{-}(M_1)=(\delta_1,-\delta_1,-1,\delta_4;\sigma_1 e_4,\sigma_2 e_3)$$
 where $\delta_i,\sigma_j \in \{\pm1\}$ are arbitrary for $i=1,4,\, j=1,2$.
Similarly, condition ($\vep_1$) for $M_1'$ gives $\delta_2=\pm1$,
$\delta_1\delta_4=\pm$ and $\delta_1\delta_2=-1$. Thus, $M_1'$ has $2^3$
$\pinpm$ structures given by
$$\vep_\pm(M_1')=(\mp1,\pm1,\delta_3,-1;\sigma_1 e_4,\sigma_2 e_3)$$
 with $\delta_3,\sigma_1,\sigma_2\in\{\pm1\}$. In this way, we have shown
  that $M_1,M_1'$ is an isospectral pair such that $M_1$ carries no $\pinp$ structure while ${M_1}'$ admits $2^3$ of them.

We note that the orientable manifolds $\tilde{M}_1,{\tilde{M}_1}'$ do have the same  properties.
These manifolds are still isospectral (again we have equality of  Sunada numbers) and $\g_i^2$ and  ${\g_i'}^2$ are the same as before, for $1\le i \le 3$.

Now, if we look for spin structures $\vep$ on $\tilde{M}_1, {\tilde{M}_1}'$, we get
\begin{gather*}
u^2(B_1)=u^2(B_1')=(\sigma_1 e_4e_5)^2= -1, \quad u^2(B_2)= u^2(B_2')= (\sigma_2 e_3e_6)^2=-1, \\ u^2(B_3)=u^2(B_3')=(\sigma_3 e_3e_4e_5 e_6)^2=1.
\end{gather*}
For $\tilde{M}_1$ we have $\g_2^2= \g_3^2 =L_{e_1+e_2}$, hence $\vep (\g_2^2)= \vep(\g_3^2)$, a contradiction, given that $u^2(B_2)= -1$ and $u^2(B_3)=1$.
Thus,  there are no spin structures on $\tilde{M}_1$.
On the other hand, for ${\tilde{M}_1}'$, we have ${\g_1'}^2=L_{e_2}, \: {\g_2'}^2=L_{e_1+e_4}, \: {\g_3'}^2=L_{e_1+e_2}$.
Thus,  $\vep (L_{e_2}) = -1 ,\: \vep (L_{e_1+e_4}) = -1,\: \vep (L_{e_1+e_2}) = 1$, hence there are
$2^5$ spin structures given by
$$\vep=(-1,-1,\delta_3,1,\delta_5, \delta_6;\sigma_1 e_4 e_5,\sigma_2 e_3 e_6)$$
 with $\delta_3,\delta_5,\delta_6,\sigma_1,\sigma_2\in\{\pm1\}$.

This proves the claim and shows that one can not hear the existence of spin structures on a compact Riemannian manifold.

\medskip
We consider next the remaining pairs $M_i, M_i'$, $2\leq i\leq 5$. The calculations 
are entirely similar to those in the cases discussed above,
so we will omit  the details, giving the necessary information in several
tables. For convenience, we will also include the pair  $M_1,M'_1$.

Note that the manifolds $M_1,M_1',M_2,M_2'$, as well as
$M_3,M_3',M_4,M_4'$, have the same holonomy representation. Furthermore,
 all matrices appearing in Table 1 are conjugate to
 $B_{0,1}$, $B_{0,2}$ or  $B_{0,3}$. By Lemma~\ref{preimages} we know
 that ${u^\pm_{0,1}}^2=\pm1$, ${u^\pm_{0,2}}^2=-1$ and
  ${u^\pm_{0,3}}^2=\mp1$ for $\Pinpm$. Thus we have:
\addtocounter{table}{1}
\begin{table}[h]
\caption{}
\begin{center}
\begin{tabular}{|c|c|c|c|}
\hline
   manifolds  & $u^2_\pm(B_1)$ &   $u^2_\pm(B_2)$ &   $u^2_\pm(B_3)$    \\ \hline
$M_1,M_1',M_2,M_2'$ & $\pm1$ & $\pm1$ & $-1$ \\ \hline
$M_3,M_3',M_4,M_4'$ & $\pm1$ & $\mp1$ & $-1$ \\ \hline
$M_5,M_5'$ & $-1$ & $-1$ & $-1$ \\ \hline
\end{tabular}
\end{center}
\end{table}

One has that $\g_i^2=L_{\ld_i}\in \Ld$. In Table 3 we give the vectors
$\ld_i$ for $1\leq i\leq 3$ and for every $M_j,M_j'$, $1\leq j \leq 5$.
\begin{table}[h]
\caption{}
\begin{center}
\begin{tabular}{|c|c|c|c|c|c|c|c|c|c|c|}
\hline
 & $M_1$ & $M_1'$ & $M_2$ & $M_2'$ & $M_3$ & $M_3'$ & $M_4$ & $M_4'$ & $M_5$ & $M_5'$  \\ \hline
$\g_1^2$ & \sz{ $e_3$} & \sz{$e_2$} & \sz{$e_3$} & \sz{$e_2$} & \sz{$e_4$} & \sz{$e_2$} & \sz{$e_1+e_2$} & \sz{$e_2+e_4$} & \sz{$e_4$} & \sz{$e_3$} \\ \hline
$\g_2^2$ & \sz{$e_1+e_2$} & \sz{$e_1+e_4$} & \sz{$e_2+e_4$} & \sz{$e_1+e_2$} & \sz{$e_4$} & \sz{$e_4$} & \sz{$e_4$} & \sz{$e_4$} & \sz{$e_4$} & \sz{$e_1$} \\ \hline
$\g_3^2$ & \sz{$e_1+e_2$} & \sz{$e_1+e_2$} & \sz{$e_2$} & \sz{$e_1$} & \sz{$e_3$} & \sz{$e_4$} & \sz{$e_4$} & \sz{$e_3$} & \sz{$e_2$} & \sz{$e_2$} \\ \hline
\end{tabular}
\end{center}
\end{table}

Using the information obtained in Tables 2 and 3 we
 get %Table 4 containing
the equations to be satisfied by the $\delta_i$'s,  resulting from
condition ($\vep_1$).

\begin{table}[h] \caption{Equations for $\delta_i$, $1\le i \le 4$.}
\begin{center}
\begin{tabular}{|c|c|c|c|}
\hline
    & $\g_1$ &   $\g_2$ &   $\g_3$    \\ \hline
$M_1$  & $\delta_3=\pm1$ & $\delta_1\delta_2=\pm1$ & $\delta_1\delta_2=-1$ \\ \hline
$M_1'$ & $\delta_2=\pm1$ & $\delta_1\delta_4=\pm1$ & $\delta_1\delta_2=-1$ \\ \hline
$M_2$  & $\delta_3=\pm1$ & $\delta_2\delta_4=\pm1$ & $\delta_2=-1$ \\ \hline
$M_2'$ & $\delta_2=\pm1$ & $\delta_1\delta_2=\pm1$ & $\delta_1= -1$ \\ \hline
$M_3$  & $\delta_4=\pm1$ & $\delta_4=\mp1$ & $\delta_3=-1$ \\ \hline
$M_3'$ & $\delta_2=\pm1$ & $\delta_4=\mp1$ & $\delta_4=-1$ \\ \hline
$M_4$  & $\delta_1\delta_2=\pm1$ & $\delta_4=\mp1$ & $\delta_4=-1$ \\ \hline
$M_4'$ & $\delta_2\delta_4=\pm1$ & $\delta_4=\mp1$ & $\delta_3=-1$ \\ \hline
$M_5$  & $\delta_4=-1$ & $\delta_4=-1$ & $\delta_2=-1$ \\ \hline
$M_5'$ & $\delta_3=-1$ & $\delta_4=-1$ & $\delta_2=-1$ \\ \hline
\end{tabular}
\end{center}
\end{table}

 By looking at Table~4 we immediately see that
$M_1$ has no $\pinp$ structures, $M_2'$ and $M_3$ admit no $\pinpm$
structures and  $M_3'$ has no $\pinm$ structures, since the corresponding
equations are not compatible.
 We now list all the characters ${\vep_\pm}_{|\Ld}$,
  corresponding to the $\pinpm$ and spin  structures in the remaining
  cases.
\begin{align*}
\vep_-(M_1)&=(\delta_1,-\delta_1,-1,\delta_4), & \qquad
\vep_{\pm}(M_1')&=(\mp1,\pm1,\delta_3,-1),   \\
\vep_{\pm}(M_2)&=(\delta_1,-1,\pm1,\mp1), & \qquad
\vep_{-}(M_3)&=(\delta_1,-1,\delta_3,-1),   \\
\vep_+(M_4)&=(\delta_1,\delta_1,\delta_3,-1), & \qquad
\vep_{\pm}(M_4')&=(\delta_1,-1,-1,\mp1),   \\
\vep(M_5)&=(\delta_1,-1,\delta_3,-1), & \qquad
\vep(M_5')&=(-1,-1,-1,\delta_4).
\end{align*}

Now, for each choice of ${\vep_\pm}_{|\Ld}$
there are $2^2=4$ structures corresponding to the possible choices of $\sigma_1, \sigma_2$, 
hence it is easy to verify that  the number of $\pinp$, $\pinm$ or spin structures
is as indicated in the theorem.
\end{proof}

\begin{rem} The procedure of adding appropriate characters to $M_1,M_1'$ to obtain  orientable manifolds, with $M_1$ admitting a spin structure while $M'_1$ does not, can be used  with the remaining pairs $M_i,M_i'$, $2\leq i \leq 4$, as well.  Alternately, we can also use the method 
described in Remark~\ref{doubling}. Indeed, consider the orientable $\Z_2^2$-manifolds $M_{\text{d}\G_1},M_{\text{d}\G'_1}$ of dimension 8
 obtained by doubling the Bieberbach groups  $\G_1,\G'_1$ (see Table 5). The resulting  manifolds now carry a K\"ahler structure.
\begin{table}[h] \caption{}
\begin{center}
\begin{tabular}{|rcc|rcc|rcc|}  \hline    $B_1$ &  $L_{b_1}$  & $L_{b_1'}$    & $B_2$ &  $L_{b_2}$ & $L_{b_2'}$  & $B_3$ &  $L_{b_3}$ & $L_{b_3'}$   \\ \hline
 1 && & 1 &$\text{{\scriptsize 1/2}}$ & $\text{{\scriptsize 1/2}}$ & 1 &$\text{{\scriptsize 1/2}}$ &$\text{{\scriptsize 1/2}}$ \\
      1 && $\text{{\scriptsize 1/2}}$ & 1 & $\text{{\scriptsize 1/2}}$ & & 1 &$\text{{\scriptsize 1/2}}$ &$\text{{\scriptsize 1/2}}$  \\
 1 &  $\text{{\scriptsize 1/2}}$ & &  -1 && & -1 &  $\text{{\scriptsize 1/2}}$ &  \\
 -1 &  & & 1 & &$\text{{\scriptsize 1/2}}$ & -1 & &$\text{{\scriptsize 1/2}}$  \\
 \hline
 1 && & 1 &$\text{{\scriptsize 1/2}}$ & $\text{{\scriptsize 1/2}}$ & 1 &$\text{{\scriptsize 1/2}}$ &$\text{{\scriptsize 1/2}}$ \\
      1 && $\text{{\scriptsize 1/2}}$ & 1 & $\text{{\scriptsize 1/2}}$ & & 1 &$\text{{\scriptsize 1/2}}$ &$\text{{\scriptsize 1/2}}$  \\
 1 &  $\text{{\scriptsize 1/2}}$ & &  -1 && & -1 &  $\text{{\scriptsize 1/2}}$ &  \\
 -1 &  & & 1 & &$\text{{\scriptsize 1/2}}$ & -1 & &$\text{{\scriptsize 1/2}}$  \\
 \hline
\end{tabular}
\end{center}
\end{table}

By comparing  the Sunada numbers, we see that $M_{\text{d}\G_1}$ and $M_{\text{d}\G'_1}$ are isospectral.
Now, we look at condition ($\vep_1$) in (\ref{conditions}). For $\text{d}\G_1$ we have that $\delta_3\delta_7=-1$, $\delta_1\delta_2\delta_5\delta_6=-1$ and $\delta_1\delta_2\delta_5\delta_6=1$. These last two equations are clearly not compatible, hence $M_{\text{d}\G_1}$ admits no spin structures. On the other hand, for $\text{d}\G'_1$ we get $\delta_2\delta_6=-1$, $\delta_1\delta_4\delta_5\delta_8=-1$ and $\delta_1\delta_2\delta_5\delta_6=1$,  hence
$\vep=(\delta_1,\delta_2,\delta_3,\delta_4,-\delta_1,-\delta_2,\delta_7,-\delta_4,\sigma_1 e_4e_8,\sigma_2 e_3e_7)$, thus obtaining $2^7$ spin structures in this case.

\end{rem}

\section{Pin$^\pm$ structures on $\Z_2$-manifolds.}
In this last section we study in some detail the special case of $\Z_2$-manifolds, where an explicit description of the $\pinpm$ structures can be given.
For each $0\leq j,h < n$, let as in (\ref{Bjh})
$$B_{j,h}:=\text{diag}(\underbrace{J,\dots,J}_j,\underbrace{-1,\dots,-1}_h,\underbrace{1,\dots,1}_l)$$
 where $n=2j+h+l$, $j+h\not=0$ and $l\geq1$.  Then $\bjh \in \on$, $\bjh^2=\I$ and $\bjh \in \son$ if and only if $j+h$ is even.
Let $\Ld=\Z e_1\oplus \cdots \oplus \Z e_n$ be the canonical lattice of $\R^n$ and for $j,h$ as before define the groups
\begin{equation}
\G_{j,h}:=\langle \bjh L_{\frac{e_n}{2}}, \Ld \rangle.
\end{equation}
 We have that $\Ld$ is stable by $B_{j,h}$ and $(\bjh +\I)\frac{e_n}2 = e_n \in \Ld \smallsetminus (\bjh +\I)\Ld$. Hence, by Proposition~2.1  in \cite{DM}, the $\G_{j,h}$ are Bieberbach groups. In this way, if $M_{j,h}= \G_{j,h} \backslash \R^n$, we have a family
\begin{equation}\label{family}
\mathcal{F}=\{ M_{j,h} \,:\, 0 \le j \le [\tfrac{n-1}2], 0 \le h < n-2j, \,j+h\not=0 \}
\end{equation}
of compact flat manifolds  with holonomy group $F\simeq \Z_2$.
The next proposition summarizes some known results on $\Z_2$-manifolds.
We include a proof for completeness.

\begin{prop}
The family $\mathcal{F}$ gives a  system of representatives for the diffeomorphism classes of $\Z_2$-manifolds of dimension $n$. Furthermore we have:
\begin{equation}\label{H1}
H_1(M_{j,h},\Z) \simeq \Z^{j+l} \oplus \Z_2^h.
\end{equation}
For $1\le p\le n$,
\begin{equation}\label{bettip}
\beta_p(M_{j,h})=  \sum _{i=0}^{[\frac p2]}\binom{j+h}{2i} \binom{j+l}{p-2i}.
\end{equation}
If $\beta_1(M_{j,h})=\beta_1(M_{j',h'})$, then $\beta_p(M_{j,h})=\beta_p(M_{j',h'})$ for any $p \ge 1$.
\end{prop}
\begin{proof}
We first prove that the manifolds $\mjh$ are pairwise non homeomorphic.
We now compute $H_1(M_{j,h},\Z) \simeq \G_{j,h}/[\G_{j,h},\G_{j,h}]$. For
$\g=\bjh L_{\frac{e_n}{2}}$, we have
\begin{eqnarray*}
[\G_{j,h},\G_{j,h}]&=& \langle [\g, L_{e_i}]
= L_{(B-\I){e_i}}: 1 \le i \le n \rangle\\ & =& \langle L_{e_2- e_1},\dots, L_{e_{2j}- e_{2j-1}}, L_{2e_{2j+1}},\dots, L_{2e_{2j+h}} \rangle.
\end{eqnarray*}
Using this information and the fact that $\g^2 =L_{e_n}$ it is  easy to see that
$$H_1(M_{j,h},\Z) \simeq \Z^{j+l} \oplus \Z_2^h.$$
Thus, if $M_{j,h}$ and $M_{j',h'}$ are homeomorphic then $h=h'$ and
$j+l=j'+l'$, hence  $j=j'$ as asserted.

To show that the family $\mathcal{F}$ gives a complete system of representatives for the diffeomorphism classes of $\Z_2$-manifolds, we will use results in \cite{Ch}, p.\@153 (it could also be proved directly by using that any integral representation of $\Z_2$ decomposes uniquely as a sum of indecomposable representations of rank $\leq 2$ given by $1, -1$ or $J$).

The cardinality of $\fa$ equals 
$ \Big(\sum_{j=0}^{[\frac{n-1}2]} n-2j \Big) -1$, since we must exclude the case $j=h=0$ corresponding to $B_{0,0}=\I$.
Thus we have
\begin{equation}\label{e.card}
\#\fa  =  
\big(n-[\tfrac{n-1}2]\big)\big([\tfrac{n-1}2]+1\big)-1
=\left\{ \begin{array}{lr}\frac{n^2+2n-4}{4} & \qquad \text{$n$ even} \\ \frac{n^2+2n-3}{4} & \qquad \text{$n$ odd.} \end{array} \right.
\end{equation}

On the other hand, if $p$ is a prime, Charlap gives a formula for the number $N_p$ of diffeomorphism classes of $\Z_p$-manifolds of dimension $n$.
For $p=2$ this number is given by:
\begin{eqnarray*}
N_2 & =  & \tfrac12 [\tfrac{n-1}2]\Big([\tfrac{n-1}2]+3\Big) +\tfrac12 \Big( (n-1)-[\tfrac{n-1}2] \Big)\Big( n-[\tfrac{n-1}2] \Big).
\end{eqnarray*}
In this way we obtain that $N_2=\tfrac18(n-2)(n+4) + \tfrac18 n(n+2) = \tfrac{n^2+2n-4}{4}$ for $n$ even, and
$N_2=\tfrac14 (n-1)(n+3)=\tfrac{n^2+2n-3}4$, for $n$ odd.
This shows that $\#\fa=N_2$, as claimed.

To determine the  $p$-Betti number of $\mjh$ for $1\leq p\leq n$,  we
note that $\bjh$ acts diagonally on the basis
 $e_1 \pm e_2, \dots, e_{2j-1} \pm e_{2j}, e_{2j+1}, \dots, e_n$,  with
  $j+l$ (resp.\@ $j+h$) eigenvectors with eigenvalue 1 (resp.\@ $-1$).
   Thus, an exterior product of $p$ elements of this basis will be
   invariant by $\bjh$, if and only if  an even number of them have
    eigenvalue $-1$.  Hence we have
\begin{equation*}
\beta_p(M_{j,h})=  \sum _{i=0}^{[\frac p2]}\binom{j+h}{2i} \binom{j+l}{p-2i}
\end{equation*}
as asserted. Now,  if $\beta_1(M_{j,h})=\beta_1(M_{j',h'})$ then
$j+l=j'+l'$ and hence $j+h=j'+h'$. Thus,
$\beta_p(M_{j,h})=\beta_p(M_{j',h'})$, for any $1\le p \le n$.
\end{proof}

The next result gives a description of $\pinpm$ structures on $\Z_2$-manifolds.

\begin{prop}
 Every $\Z_2$-manifold $M_\G$ has  $\pinpm$
structures (and  spin structures, if $M_\G$ is orientable).
If $\,\G = \G_{j,h}$ then $M_\G$ has $2^{n-j}$ $\pinpm$ structures
parametrized by the tuples $(\delta_1, \dots,\delta_n, \sigma) \in \{\pm
1\}^{n+1}$
satisfying:
\begin{equation}\label{relations}
\delta_1=\delta_2, \;\cdots, \; \delta_{2j-1}=\delta_{2j}
\end{equation} and
\begin{equation}  \label{prodcond1} \delta_n=\left\{ \begin{array}{ll}
(-1)^{jh}(-1)^{[\frac{j}2]}(-1)^{[\frac{h}2]} & \quad \text{for $\pinp$ structures} \\
(-1)^{jh}(-1)^{[\frac{j+1}2]}(-1)^{[\frac{h+1}2]} & \quad \text{for $\pinm$ structures.}
\end{array} \right.
\end{equation}
In particular, in the case of spin structures we have $\delta_n=(-1)^{\frac{j+h}2}$.
\end{prop}
\begin{proof}
In  light of Proposition 4.1, we have that $\G \simeq \G_{j,h}$ for some $j,h$, hence $M_\G$ is diffeomorphic to $M_{j,h}$. Therefore, since $\pinpm$ structures on diffeomorphic manifolds are in a bijective correspondence, we may assume that $\G = \G_{j,h}$.

We have observed in Remark 2.2 that equation ($\vep_2$)  always holds for $\Z_2$-manifolds of diagonal type. However, in the non diagonal case,  ($\vep_2$) gives a restriction. Namely, let $\ld=\sum_{i=1}^n m_i e_i$, $m_i\in \Z$. Then
$$(B_{j,h}-\I)\ld=\sum_{i=1}^j(m_{2i}-m_{2i-1}) e_{2i-1}+(m_{2i-1}-m_{2i}) e_{2i} -2\sum_{i=1}^h m_{2j+i} e_{2j+i}.$$
Thus, ($\vep_2$) holds if and only if
\begin{equation*}
\delta_1^{(m_2-m_1)}\delta_2^{(m_1-m_2)}\cdots \delta_{2j-1}^{(m_{2j}-m_{2j-1})}\delta_{2j}^{(m_{2j-1}-m_{2j})}=1
\end{equation*}
for every $m_1,\dots,m_{2j}\in \Z$, or equivalently,
 \begin{equation*}\label{pairs}
 \delta_1=\delta_2,\;\ldots,\;\delta_{2j-1}=\delta_{2j}.
\end{equation*}
Each of these relations divides by 2 the number of structures. Hence we
 obtain a maximum of $2^{n-j+1}$ $\pinpm$ structures for $M_{j,h}$.
Furthermore, equation ($\vep_1$)  gives another restriction  since
$\vep_\pm(\g^2)=\vep_\pm(L_{(B+\I)b})= \vep_\pm(\g)^2$.
Now $(B+\I)b=e_n$, hence, by (\ref{squares}),
 equation ($\vep_1$) reads:
\begin{equation} \label{cond1}
\delta_{n} = \left\{ \begin{array}{ll}
(-1)^{jh}(-1)^{[\frac{j}2]}(-1)^{[\frac{h}2]} & \quad \text{in $Cl^+(n)$} \\
(-1)^{jh}(-1)^{[\frac{j+1}2]}(-1)^{[\frac{h+1}2]} & \quad \text{in $Cl^-(n).$}
\end{array} \right.
\end{equation}
Thus, the restriction imposed by (\ref{cond1}) divides by 2 the number of structures and we get a total of $2^{n-j}$ pin$^\pm$ structures on $\man$ for $\G= \G_{j,h}$. 
\end{proof}

\noindent {\em Note.}
 Proposition 4.1 together with Lemma 3.1, give an
explicit description of all $\pinpm$ structures on  $\Z_2$-manifolds.

\begin{ejem}
As a final task, to illustrate Proposition 4.2, we  list explicitly the 28
$\pinpm$ Riemannian $\Z_2$-manifolds $(M,\vep)$ of dimension 3 having
canonical lattice of translations $\Ld$.

There are  3 diffeomorphism classes,  one of which splits into 2 isometry classes, hence we have 4 isometry classes, corresponding to the groups 
  $\G_{1,0}= \left \langle \left[\begin{smallmatrix} J  & \\   & 1 \end{smallmatrix}\right]L_{\frac{e_3}2},\, \Ld \right\rangle$,
$\G_{0,1}=\left\langle \left[\begin{smallmatrix} -1 &&\\&1&\\&&1 \end{smallmatrix}\right]L_{\frac{e_3}2},\, \Ld \right\rangle$,
$\G_{0,1}'= \left\langle \left[\begin{smallmatrix} -1 &&\\&1&\\&&1 \end{smallmatrix}\right]L_{\frac{{e_2}+{e_3}}2},\, \Ld \right\rangle$ and  $\G_{0,2}=\left\langle \left[\begin{smallmatrix} -1 &&\\&-1&\\&&1 \end{smallmatrix}\right]L_{\frac{e_3}2},\,  \Ld \right\rangle$.

We note that $M_{0,1}$ and $M_{0,1}'$ are not isometric, as can be seen by computing the injectivity radius, that is the length of the shortest closed geodesic. Indeed, using the results in \cite{MR3} one easily sees that these equal $\frac 12$ and $\frac {\sqrt 2}2$, respectively.

The pairs $(M,\vep)$ are listed in the following table, obtained by using Lemma \ref{preimages} and Theorem \ref{z2spinst}.

\begin{center} {\sc Table 6.} Pin$^\pm$ structures on $\Z_2$-manifolds of dimension 3.
 \smallskip

\begin{tabular}{|c|c|c|l|c|}
\hline
$M_\G$ & cond.\@ ($\vep_1$) & cond.\@ ($\vep_2$) & $\pinpm$ structures & $\#$\\
 \hline
$M_{1,0}$ & $\delta_3=\pm1$ & $\delta_1=\delta_2$ & $\vep_\pm=(\delta_1,\delta_1,\pm1;\sigma\frac{\sqrt 2}{2}(e_1 - e_2))$ & $2^2$ \\ \hline
$M_{0,1}$ & $\delta_3=\pm1$ & $-$ & $\vep_\pm=(\delta_1,\delta_2,\pm1;\sigma e_1)$ & $2^3$ \\ \hline
$M_{0,1}'$ & $\delta_2\delta_3=\pm1$ & $-$ & $\vep_\pm=(\delta_1,\delta_2,\pm\delta_2;\sigma e_1)$ & $2^3$ \\ \hline
$M_{0,2}$ & $\delta_3=-1$ & $-$ & $\vep_\pm=(\delta_1,\delta_2,-1;\sigma e_1e_2)$ & $2^3$ \\ \hline
\end{tabular}
\end{center}

\medskip
We note that the spin structures for $M_{0,2}$  are already contained in \cite{Pf}.

\end{ejem}

                % End The Bibliography

\end{document}